\documentclass[11pt,intlimits, draft,oneside]{amsart}

\usepackage[latin2]{inputenc}
\usepackage{latexsym}
\usepackage{amssymb}
\usepackage{amsthm}

\usepackage{color} % Farbe

\usepackage{amsmath}

\usepackage{enumerate}

\input{cyracc.def}
\newfont{\cyrfnt}{wncyi10 at 11pt}

\newtheorem{thm}{Theorem}[section]
\newtheorem{prop}[thm]{Proposition}
\newtheorem{lemma}[thm]{Lemma}

\theoremstyle{definition}

\newtheorem{remark}[thm]{Remark}

\newcommand{\R}{\mathbb{R}}

\newcommand{\N}{\mathbb{N}}
\newcommand{\Z}{\mathbb{Z}}
\newcommand{\U}{\mathcal{U}}
\begin{document}

\title[Category theorems for stable semigroups]
{Category theorems for stable semigroups}
%\title[Category theorems for stable operators]{Category theorems for weakly and almost weakly stable operators on Hilbert spaces}
%\footnote{The fourth author was supported by the Marie Curie Host Fellowship "Spectral theory for evolution equations", contract number HPMT-CT-2001-00315.}

\author{Tanja Eisner}
\address{Tanja Eisner \newline Mathematisches Institut, Universit\"{a}t T\"{u}bingen\newline Auf der Morgenstelle 10, D-72076, T\"{u}bingen, Germany}
\email{talo@fa.uni-tuebingen.de}

\renewcommand{\thefootnote}{}
\author{Andr\'{a}s Ser\'{e}ny}
\footnote{The second author was supported by the DAAD-PPP-Hungary Grant, project number D/05/01422.}
%contract number HPMT-CT-2001-00315.}
\address{Andr\'{a}s Ser\'{e}ny \newline Department of Mathematics and its Applications, Central European University,\newline N\'{a}dor utca  9, H-1051 Budapest, Hungary }
\email{sandris@elte.hu}

\keywords{Weak and strong mixing, category, $C_0$-semigroups, Hilbert space, stability}
%\subjclass[2000]{47A35, 37A25}

\begin{abstract}
Inspired by the classical category theorems of Halmos and Rohlin for
discrete measure preserving transformations, we prove analogous
results in
%extend these results to
the abstract setting of unitary and isometric $C_0$-semigroups on a
separable Hilbert space. More presicely, we show that the set of all weakly stable
unitary groups (isometric semigroups) is of first category, while the set of all almost weakly stable
unitary groups (isometric semigroups) is residual for an appropriate topology.
%
%We compare the two related notions of weak and almost weak stability for
%unitary and isometric $C_0$-(semi)groups on a Hilbert space and
%show that, in the Baire category sense, these notions differ
%fundamentally. More precisely, we show that the set of all weakly stable
%unitary groups (isometric semigroups) is of first and the set of all almost weakly stable
%unitary groups (isometric semigroups) is residual for an appropriate
%topology. This gives an extension of classical category theorems of
%Halmos and Rohlin for measure preserving transformations in ergodic
%theory.
\end{abstract}

\maketitle

\vspace{0.1cm}

%%%%%%%%%%%%%%%%%%%%%%%%%%%%%%%%%%%%%%%%%%%%%%%%%%%%%%%%%%%%%%%%%%%%%%%
%%                                                                   %%
%%             INTRODUCTION                                          %%
%%                                                                   %%
%%%%%%%%%%%%%%%%%%%%%%%%%%%%%%%%%%%%%%%%%%%%%%%%%%%%%%%%%%%%%%%%%%%%%%%

\section{Introduction}

In 1944 Halmos \cite{halmos:1944} showed that a ``typical''
dynamical system is weakly mixing, or, more precisely, the set of
all weakly mixing transformations on a measure space is residual and
hence of second category. Four years later, Rohlin
\cite{rohlin:1948} showed that the set of all strongly mixing
transformation is of first category (see also Halmos
\cite[pp.~77--80]{halmos:1956}). In this way they proved the
existence of weakly but not strongly mixing systems (see e.g. Petersen
\cite[Section 4.5]{petersen:1983} for later concrete examples and
for a method to construct such systems). 
%%%Their result has been reproduced many times (see ... //B\"ucher!//). 
%%%All these authors consider the discrete time case only, while it is common belief that a similar
%%%result should hold for  continuous time dynamical systems as well.

%There are several analogues of the category theorems for the time continuous case. For example, recently, Bartoszek and Kuna \cite{bartoszek/kuna:2006} proved an analogue of the Halmos and Rohlin result 
In the time continuous case, i.e., for measure preserving flows, analogous results %to hold but are 
%seem not to be 
are not stated in the standard literature. However, Bartoszek and Kuna \cite{bartoszek/kuna:2006} recently proved a category theorem for Markov semigroups on the Schatten class $C_1$, while an analogous result for stochastic semigroups can be found in Lasota, Myjak \cite{lasota/myrjak:1992}.
%on stochastic semigroups
%Iwanik \cite{iwanik:1992}...

In this paper we consider the time continuous case in the more general setting of unitary and isometric semigroups on a Hilbert space and prove corresponding category results. 

%In this paper we treat the time continuous case in a more general
%setting and make precise in which sense and for which topology such
%a category result holds.

Recall that a (continuous) dynamical system $(\Omega, \mu, (\varphi_t)_{t\geq 0})$ induces a
unitary/ isometric $C_0$-(semi)group $(T(t))_{t\geq 0}$ on the Hilbert
space $H:=L^2(\Omega,\mu)$ via the formula $T(t)f(\omega)=f(\varphi_t (\omega))$. Moreover, the decomposition into invariant subspaces $H=\langle \mathbf{1}\rangle \oplus H_0$ holds for $H_0=\{f\in H: \int_\Omega f d\mu=0\}$. 
Recall further that the flow $(\varphi_t)$ is strongly mixing if and
only if the restricted semigroup $(T_0(t))_{t\geq 0}$ on $H_0$ is weakly stable, i.e.,
$\lim_{t\to\infty}\langle T_0(t)x,y\rangle=0$ for every $x,y\in H_0$. Analogously, weak mixing of the flow corresponds to
so-called \emph{almost weak stability} of the semigroup $T_0(\cdot)$,
i.e.,
\begin{equation*}
\lim\langle T_0(t)x,y\rangle=0, \quad t\to\infty,\ t\in M, \quad \text{for every } x,y\in H_0
\end{equation*}
for a set $M\subset \R_+$ with (asymptotical) density one, i.e., with
$d(M):=\lim_{t\to\infty}\frac{\lambda(M\cap[0,t])}{t}=1$ for the Lebesgue
measure $\lambda$. We refer to e.g. Cornfeld, Fomin, Sinai \cite[Section 1.7]{cornfeld/fomin/sinai:1982} for details.

We now consider an arbitrary Hilbert space $H$ and a contractive
$C_0$-semigroup $T(\cdot)$ on $H$ with generator $A$. By the
classical decomposition theorem of Jacobs--Glicksberg--de Leeuw (see
Glicksberg, De Leeuw \cite{glicksberg/deleeuw:1961} or Engel, Nagel
\cite[Theorem V.2.8]{engel/nagel:2000}), $H$ decomposes into the
so-called reversible and stable parts
\begin{eqnarray*}
H_r&:=&\overline{\text{lin}}\{x: Ax=i\alpha x \text{ for some } \alpha
\in\R\}, \\
H_s&:=&\{x: 0 \text { is a weak accumulation point of } \{T(t)x,t\geq 0\}\}.
\end{eqnarray*}
(Note that the Jacobs--Glicksberg--de Leeuw decomposition is valid for
every bounded semigroup on a reflexive Banach space or, more generally,
for every relatively weakly compact semigroup.)
A result of Hiai \cite{hiai:1978} being a continuous version of
discrete results of Nagel \cite{nagel:1974} and Jones, Lin
\cite{jones/lin:1976} shows that the semigroup $T(\cdot)$ restricted
to $X_s$ is almost weakly stable as defined above. (For a survey
on weak and almost weak stability see Eisner, Farkas, Nagel, Serény \cite{EFNS}.)

%So we see that almost weak stability is an important asymptotic
%property of $C_0$-semigroups appearing in a very general setting and
%being equivalent to emptiness of the point spectrum of the generator on the
%imaginary axis.
%On the other side, weak stability also being an important stability
%property appearing in other decomposition theorems, see Foguel
%\cite{foguel:1963}, is much more difficult to treat. This is
%still an open question how one can characterise weak stability of a
%$C_0$-semigroup in terms of its generator. Even for unitary groups
%(which are never strongly stable if non-constant) with
%bounded generators this question does not become easier.  (Quite
%recently, Chill and Tomilov \cite{chill/tomilov:2004} introduced a
%sufficient condition, see also Eisner, Farkas, Nagel, Serény
%\cite{EFNS}, but it is not clear for which semigroups such
%conditions are also necessary.) Note that this difficulty appears in
%ergodic theory as well and corresponds to investigation of strong mixing.

Weak mixing (almost weak stability) can be characterised by a simple
spectral condition, while the more natural property of strong mixing
(weak stability) \emph{``is, however, one of those notions, that is
easy and natural to define but very difficult to study...''}, see
the monograph of Katok, Hasselblatt
\cite[p.~748]{katok/hasselblatt:1995}.

In this paper we extend the results of Halmos and Rohlin and show
that weak and almost weak stability differ fundamentally. More
precisely, we prove that for an appropriate (and natural) topology
the set of all weakly stable unitary groups on a separable
infinite-dimensional Hilbert space is of the first category and the
set of all almost weakly stable unitary groups is residual and hence
of the second category (Section 2). An analogous result holds for
isometric semigroups (Section 3). In Section 4 we discuss the case
of contractive semigroups.

For analogous results on the time discrete case, see Eisner,
Ser\'eny \cite{eisner/sereny:2006}.

\section{Unitary  case} \label{section:unitary-gr}

%Let $H$ be a separable infinite-dimensional Hilbert space.
%We denote the set of all unitary operators on $H$ by $\U$.
%The following density result for periodic operators is a first
%building block for our construction.

Consider the set $\U$ of all unitary $C_0$-groups on an infinite-dimensional Hilbert
space $H$. The following density result on periodic unitary groups is a first step in our construction.

\begin{prop}\label{prop:periodic} For every $N\in\N$
and every unitary group $U(\cdot)$ on $H$ there is a sequence $\{V_n(\cdot)\}_{n=1}^\infty$ of periodic unitary groups with period greater than $N$ such that $\lim_{n\to\infty}\|U(t)-V_n(t)\|=0$ uniformly on compact time intervals.
%the set of all
%  periodic unitary groups with period greater than $n$ is dense in
%  $\U$ endowed with the operator norm topology uniform on compact time intervals.
\end{prop}
\begin{proof}

Take $U(\cdot)\in \U$ and $n\in \N$.
%and $\varepsilon>0$.
By the spectral
theorem (see, e.g., Halmos \cite{halmos:1963}), $H$ is isomorphic to $L^2(\Omega, \mu)$ for some locally
compact space $\Omega$ and measure $\mu$ and $U(\cdot)$ is
unitary equivalent to a multiplication group $\tilde{U}(\cdot)$ with
\begin{equation*}
(\tilde{U}(t)f)(\omega)=e^{itq(\omega)}f(\omega),\quad \forall
\omega\in\Omega,\ t\in \R, \ f\in L^2(\Omega, \mu)
\end{equation*}
for some measurable $q:\Omega \to \R$.

We approximate the (unitary) group $\tilde{U}(\cdot)$ as follows. Take $n>N$ and
define
\begin{equation*}
q_n(\omega):= \frac{2\pi j}{n},\quad \forall \omega\in q^{-1}\left(\left[\frac{2\pi
  j}{n}, \frac{2\pi (j+1)}{n}\right]\right),\ j\in\Z.
\end{equation*}
Denote now by $\tilde{V}_n(t)$ the multiplication operator with
$e^{itq_n(\cdot)}$ for every $t\in\R$.
The unitary group $\tilde{V}_n(\cdot)$ is periodic with period greater
than or equal to $n$ and therefore $N$. Moreover,
\begin{eqnarray*}
\Vert \tilde{U}(t)-\tilde{V}_n(t) \Vert \leq \sup_{\omega} |e^{itq(\omega)} - e^{itq_n(\omega)}| 
\leq |t| \sup_{\omega} |q(\omega)-q_n(\omega)| \leq \frac{2\pi |t|}{n} \underset{n\to\infty}{\longrightarrow} 0
%\Vert \tilde{U}(t)f - \tilde{V}_n(t)f \Vert &=&
%\left(\int_\Omega |e^{itq(\omega)} - e^{itq_n(\omega)}|^2 \|f(\omega)\|^2
%d\omega\right)^{1/2} \\
%&\leq& 
%|t| \sup_{\omega} |q(\omega)-q_n(\omega)|\cdot \|f\| \leq \frac{2\pi |t|}{n} \|f\|
\end{eqnarray*}
%holds. So $\lim_{n\to\infty}\|\tilde{U}(t) - \tilde{V}_n(t)\|= 0$
uniformly in $t$ on compact intervals and the proposition is proved.
\end{proof}
\begin{remark}
By a modification of the proof of Proposition \ref{prop:periodic}
one can show that for every $N\in\N$ the set of all
  periodic unitary groups with period greater than $N$ with bounded generators is dense in
  $\U$ with respect to the strong operator topology uniform on compact
  time intervals.
\end{remark}
%
%Before we present a second building block we need the following lemma.

For the second step we need the following lemma. From now on we assume the Hilbert space $H$ to be separable.
\begin{lemma} \label{lemma:discrete-appr-of-I}
Let $H$ be a separable infinite-dimensional Hilbert space. Then there
exists a sequence $\{U_n(\cdot )\}_{n=1}^\infty$ of almost
weakly stable unitary groups with bounded generator satisfying
$\lim_{n\to\infty}\Vert U_n(t) - I\Vert= 0$ uniformly in $t$ on compact intervals.
\end{lemma}

\begin{proof}
By the isomorphism of all separable infinite-dimensional Hilbert spaces
we can assume without loss of generality that $H=L^2(\R)$ with respect to the Lebesgue measure.
%there exists a unitary operator $U:H\to L^2(\R)$, where $L^2(\R)$ is considered with the Lebesgue measure.

Take $n\in\N$ and define $U_n(\cdot)$ on $L^2(\R)$ by
\begin{equation*}
(U_n(t)f)(s):= e^{\frac{itq(s)}{n}}f(s), \ \ s\in\R,\ \ f\in L^2(\R),
\end{equation*}
where $q:\R\to (0,1)$ is a strictly monotone increasing function.

Then all $U_n(\cdot)$ are almost weakly stable by the theorem of Jacobs--Glicksberg--de Leeuw and we have
\begin{equation*}
\Vert U_n(t) - I \Vert = \sup_{s\in\R} |e^{\frac{itq(s)}{n}} -
1| \leq [\text{for } t\leq \pi n]\leq |e^{\frac{it}{n}} - 1| \leq \frac{2t}{n} \to 0, \quad n\to \infty,
\end{equation*}
uniformly on t in compact intervals.
%
%To finish the proof we only need to define $T_n(t):=U^*\tilde{T}_n(t)U$ on $H$.
\end{proof}

%?? We now introduce the appropriate topology.
% We say that a sequence
%$\{T_n\}\subset \mathcal{L}(H)$ converges to $T\in \mathcal{L}(H)$ in
%the strong*-topology if $T_n\to T$ and $T_n^*\to T^*$ strongly (for
%details see, e.g., Takesaki \cite[p.~68]{takesaki}.
%Further we consider the space $\mathcal{U}$ of all unitary operators on $H$ endowed with the strong*-topology. Note that $\U$ is a complete metric space with respect to the metric given by
%

%Up to now we assume the Hilbert space $H$ to be separable.

The metric we introduce now on the space $\U$ is given by
%from the seminorms $p_{x,t}:=\sup_{s\in [-t,t]}\|U(s)x\|$.
%This topology is a continuous analogon of the so-called
%strong* operator topology for operators, see, e.g., Takesaki \cite[p.~68]{takesaki}). Note that the convergence for $t\to\infty$ corresponding to this topology is strong convergence uniform
%on compact time intervals of a (semi)group and its adjoint.
%from the metric
%
\begin{equation*}
d(U(\cdot),V(\cdot)):= \sum_{n,j=1}^\infty
\frac{\sup_{t\in[-n,n]}\|U(t)x_j -V(t)x_j\|}
{2^{j+n} \|x_j\|}\quad  \text{for } U,V\in \U,
\end{equation*}
where $\{x_j\}_{j=1}^\infty$ some fixed dense subset of $H$ with all
$x_j\neq 0$.
%
%This topology is a continuous analogon of the so-called
%strong* operator topology for operators, see, e.g., Takesaki \cite[p.~68]{takesaki}). Note that the convergence for $t\to\infty$ corresponding to this topology is strong convergence uniform
%on compact time intervals of a (semi)group and its adjoint.
Note that the topology coming from this metric is a continuous analogue of the so-called
strong* operator topology for operators, see, e.g., Takesaki \cite[p.~68]{takesaki}), and the corresponding convergence is
the strong convergence uniform on compact time intervals of (semi)groups and their adjoints. This metric makes $\U$ a complete metric space.
%
%We note that this metric is by the theorem of Trotter-Kato equivalent
%to the following one
%
%\begin{equation*}
%d_1(U(\cdot),V(\cdot)):= \sum_{j=1}^\infty
%\frac{\|R(1,A)x_j -R(1,B)x_j\|+\|R(1,A)^*x_j -R(1,B)^*x_j\|}
%{2^j \|x_j\|}\quad  \text{for } U,V\in \U,
%\end{equation*}
%where $A$ and $B$ are the generators of $U(\cdot)$ and $V(\cdot)$, respectively.

We denote by $\mathcal{S_U}$ the set of all weakly stable unitary groups on $H$ and by  $\mathcal{W_U}$ the set of all almost weakly stable unitary groups on $H$.
The following proposition shows the density of $\mathcal{W_U}$ which will
play an important role later.

%? Further we denote by $\mathcal{S_U}$ the set of all weakly stable unitary operators on $H$ and by $\mathcal{W_U}$ the set all almost weakly stable unitary operators on $H$.

%? We now show the following density property for $\mathcal{W_U}$.

\begin{prop}\label{prop:almweakstab}
The set $\mathcal{W_U}$ of all almost weakly stable unitary groups
with bounded generators is dense in $\U$.
\end{prop}

\begin{proof}
By Proposition \ref{prop:periodic} it is enough to approximate
periodic unitary groups by almost weakly stable unitary groups. Let $U(\cdot)$ be a
periodic unitary group with generator $A$ and period $\tau$. Take
$\varepsilon>0$, $n\in\N$, $x_1,\ldots, x_n\in H$
%\setminus\{0\}
and $t_0>0$. We have to find an almost weakly stable unitary group
$T(\cdot)$ with $\Vert U(t)x_j - T(t)x_j \Vert \leq \varepsilon$
%and $\Vert U^*x_j - T^*x_j \Vert \leq\varepsilon$
for all $j=1,\ldots, n$ and $t\in[-t_0,t_0]$.

By Engel, Nagel \cite[Theorem IV.2.26]{engel/nagel:2000} we have the orthogonal space decomposition
\begin{equation}\label{eq:periodic-decomp}
\displaystyle H=\overline{\bigoplus_{k\in\Z}
%^\bot
\ker\left(A-\frac{2\pi
      ik}{\tau}\right)}.
  \end{equation}
%where $A$ denotes the generator of $U(\cdot)$.
%
%So we can assume without loss of generality that $\{x_j\}_{j=1}^n$ is
%an orthonormal system of eigenvalues of $A$.

Assume first that  $\{x_j\}_{j=1}^n$ is an orthonormal system of eigenvectors of $A$.

Our aim is to use  Lemma \ref{lemma:discrete-appr-of-I}. For this
purpose we first construct a unitary group $V(\cdot)$ which
coincides with $U(\cdot)$ on every $x_j$ and has infinite-dimensional
eigenspaces only. %After we have done it we can apply Lemma \ref{lemma:discrete-appr-of-I}.

Define the $U(\cdot)$-invariant subspace $H_0:=\text{lin}\{x_1, \ldots,
x_n\}$ and the unitary group $V_0(\cdot):=U(\cdot)|_{H_0}$ on $H_0$. Since $H$ is separable, we can decompose $H$ in an orthogonal sum
\begin{equation*}
\displaystyle H=\bigoplus_{k=0}^\infty H_k \quad \text{with } \dim H_k=\dim H_0 \text{ for every  } k\in\N.
  \end{equation*}
%holds, where $\dim H_k=\dim H_0$ for every $k\in\N$.
%
Denote by $P_k$ an isomorphism from $H_k$ to $H_0$ and define
$V_k(\cdot):=P_k^{-1}V_0(\cdot)P_k$ on $H_k$ being copies of $V_0(\cdot)$. Consider now the unitary
group $\displaystyle V(\cdot):=\bigoplus_{k=0}^\infty V_k (\cdot)$ on $H$
which is periodic by Engel, Nagel \cite[Theorem
IV.2.26]{engel/nagel:2000}. (Note that its generator is even bounded.)

For the periodic unitary group $V(\cdot)$ the space decomposition analogous to
(\ref{eq:periodic-decomp}) holds. Moreover there are only
finitely many eigenspaces (less or equal to $n$, depending on $x_1,...,x_n$) and they all are infinite-dimensional by the
construction. Applying Lemma \ref{lemma:discrete-appr-of-I} to each
eigenspace we find an almost weakly stable unitary group $T(\cdot)$
with $\Vert U(t)x_j - T(t)x_j \Vert \leq \varepsilon$
for all $j=1,\ldots, n$ and $t\in[-t_0,t_0]$.

%From the construction follows that
%%
%\begin{equation*}
%H=\ker\left(B-\frac{2\pi i \lambda_1}{\tau}\right)\oplus^\bot \ldots
%\oplus^\bot \ker\left(B-\frac{2\pi i \lambda_n}{\tau}\right),
%  \end{equation*}
%%
%where $\frac{2\pi i \lambda_j}{\tau}$ is the corresponding eigenvalue
%of $A$ (and therefore of $B$) to the eigenvector $x_j$.
%
%Denote $X_j:=\ker\left(B-\frac{2\pi i \lambda_1}{\tau}\right)$ for
%every $j=1,\ldots,n$.
%%
%On every $X_j$
%the operator $B$ is equal to $\frac{2\pi i \lambda_j}{\tau}I$.
%Note further that all $X_j$ are infinite-dimensional.
%By Lemma \ref{lemma:discrete-appr-of-I} for every $j$ there exists an
%almost weakly stable unitary group $T_j(\cdot)$ on $X_j$ such that
%$\Vert T_j(t)-e^{\frac{2\pi t i \lambda_j}{\tau}} I\Vert <\varepsilon$
%for every $t$ with $|t|\leq t_0$. Denote now by $T(\cdot)$ the orthogonal
%sum of $T_j(\cdot)$ which is a weakly stable unitary group with
%bounded generator. Moreover, we obtain in particular that
%
%\begin{equation*}
%\Vert T(t)x_j - U(t)x_j \Vert
%= \Vert T(t)x_j - e^{\frac{2\pi i t \lambda_j}{\tau}}x_j \Vert
%\leq \varepsilon
%  \end{equation*}
%for every $t$ with $|t|\leq t_0$. % and the proposition is proved.
%

Take now arbitrary $x_1,\ldots,x_n\in H$. Take further an orthonormal basis of eigenvectors $\{y_k\}_{k=1}^\infty$. Then there exists $N\in \N$ such that
$x_j=\sum_{k=1}^N a_{jk} y_k + o_j$ with $\|o_j\|<\frac{\varepsilon}{4}$ for every $j=1,\ldots,n$.

%TO DO: korrigieren (kontinuierlich!)

We can apply the arguments above to $y_1,\ldots,y_N$ and find an almost weakly stable unitary group $T(\cdot)$ with $\|U(t)y_k-T(t)y_k\|< \frac{\varepsilon}{4N M}$
%and $\|U^*y_k-T^*y_k\|< \frac{\varepsilon}{4N M}$
for $M:=\max_{ k=1,\ldots,N, j=1,\ldots,n}{|a_{jk}|}$ and every $k=1,\ldots,N$ and $t\in[-t_0,t_0]$.
Hence
$$\|U(t)x_j-T(t)x_j\|\leq \sum_{k=1}^N |a_{jk}|\|U(t)y_k-T(t)y_k\|+2\|o_j\|< \varepsilon$$
for every $j=1,\ldots,n$ and $t\in[-t_0,t_0]$, and the proposition is proved.
\end{proof}

We are now ready to prove a category theorem for weakly and almost
weakly unitary groups.
%? We can now prove the following category theorem for weakly and almost weakly stable unitary operators.
%To do so
For this purpose we extend the argument used in the proof of the
category theorems for operators induced by measure preserving transformation in ergodic theory (see Halmos \cite[pp.~77--80]{halmos:1956}).

\begin{thm}\label{thm:unitary}
The set $\mathcal{S_U}$ of weakly stable unitary groups is of first
category and the set  $\mathcal{W_U}$ of almost weakly stable unitary groups is residual in $\U$.
\end{thm}

\begin{proof} %? First we prove that $\mathcal{S}$ is of first
              %category in $\U$.
We first prove the first part of the theorem.
 Fix $x \in H$ with $\|x\|=1$ and consider the sets

\begin{equation*}
M_t:=\left\{U(\cdot)\in \U :\ |\langle U(t) x,x\rangle| \leq \frac{1}{2} \right\}.
\end{equation*}
Note that all sets $M_t$ are closed.

For every weakly stable $U(\cdot)\in \U$ there exists $t>0$ such that
$U\in M_s$ for all $s\geq t$, i.e., $\displaystyle U(\cdot)\in N_t:=\cap_{s\geq t} M_t$. So we obtain
\begin{equation}
\mathcal{S_U} \subset \bigcup_{t>0} N_t.
\end{equation}
%where $N_t:=\cap_{s\geq t} M_s$.
Since all $N_t$ are closed, it remains to show that $\U\setminus N_t$ is dense for every $t$.

Fix $t>0$ and
let $U(\cdot)$ be a periodic unitary group. Then $U(\cdot)\notin M_s$ for some $s\geq t$ and therefore $U(\cdot)\notin N_t$. Since by Proposition \ref{prop:periodic} periodic unitary groups are dense in $\U$, $\mathcal{S}$ is of first category.

To show that $\mathcal{W_U}$ is residual  we take a dense subspace $D=\{x_j\}_{j=1}^\infty$ of $H$ and define
\begin{equation*}
W_{jkt}:=\left\{U(\cdot)\in \U :\ |\left<U(t) x_j,x_j\right>| < \frac{1}{k} \right\}.
\end{equation*}
All these sets are open. Therefore the sets
$\displaystyle W_{jk}:=\bigcup_{t>0} W_{jkt}$
are open as well.

We show that
\begin{equation}\label{W}
\mathcal{W_U}=\bigcap_{j,k=1}^\infty W_{jk}
\end{equation}
holds.

The inclusion ``$\subset$'' follows
from the definition of almost weak stability.
To prove the converse inclusion we take $U(\cdot)\notin \mathcal{W_U}$.
Then there exists $x\in H$ with $\|x\|=1$ and $\varphi \in
\R$ such that $U(t)x=e^{it\varphi} x$ for all $t>0$.
%Therefore $|\left<U(t)x,x\right>|=1$ holds for every $t\in \R$.
Take $x_j\in D$ with $\|x_j-x\|\leq \frac{1}{4}$. Then
\begin{eqnarray*}
|\left<U(t) x_j,x_j \right>| &=& |\left<U(t) (x-x_j),x-x_j \right> + \left<U(t) x,x \right>
- \left<U(t) x,x-x_j \right> - \left<U(t) (x-x_j),x \right>| \\
&\geq& 1-\|x-x_j\|^2 -2\|x-x_j\|>\frac{1}{3}
\end{eqnarray*}
for every $t>0$.
So
$U(\cdot)\notin W_{j3}$ which implies
$U(\cdot)\notin \cap_{j,k=1}^\infty W_{jk}$. Therefore equality (\ref{W}) holds.
Combining this with Proposition \ref{prop:almweakstab} we obtain that $\mathcal{W_U}$ is residual as a dense countable intersection of open sets.
\end{proof}

%%%%%%%%%%%%%%%%%%%%%%%%%%%%%%%%%%%%%%%%%%%%%%%%%%%%%%%%%%%%%%%%%%%%%%%
%%                                                                   %%
%%             ISOMETRIC CASE                                        %%
%%                                                                   %%
%%%%%%%%%%%%%%%%%%%%%%%%%%%%%%%%%%%%%%%%%%%%%%%%%%%%%%%%%%%%%%%%%%%%%%%
\section{Isometric case} \label{section:isometric-sgr}

In this section we extend the result from the previous section to isometric semigroups on $H$ which we again assume to be infinite-dimensional and separable.

Denote by $\mathcal{I}$ the set of all isometric $C_0$-semigroups on $H$. On $\mathcal{I}$ we consider the metric given by the formula
%the topology coming from the seminorms $p_{x,t}(T(\cdot)):=\sup_{s\in [0,t]} \|T(s)x\|$. The corresponding convergence is the strong convergence uniform on compact time intervals.
%Note that $\mathcal{I}$ is a complete metric space with respect to the metric given by the formula
\begin{equation*}
d(T(\cdot),S(\cdot)):= \sum_{n,j=1}^\infty \frac{\sup_{t\in[0,n]}\|T(t)x_j -S(t)x_j\|}{2^{j+n} \|x_j\|}\quad \text{for } T(\cdot),S(\cdot) \in \mathcal{I},
\end{equation*}
where $\{x_j\}_{j=1}^\infty$ is a fixed dense subset of $H$. This
corresponds to the strong convergence uniform on compact time intervals. The space $\mathcal{I}$ endowed with this metric is a complete metric space.

Analogously to Section \ref{section:unitary-gr} we denote by $\mathcal{S_I}$ the set of all weakly stable
%isometric semigroups on $H$
and by $\mathcal{W_I}$ the set all almost weakly stable isometric semigroups on $H$.

The key for our results in this section is the following classical structure
theorem on isometric semigroups on Hilbert spaces.
\begin{thm}\label{thm:Wold}({\it Wold decomposition}, see \cite[Theorem III.9.3]{sznagy/foias}.)
Let $V(\cdot)$ be an isometric semigroup on a Hilbert space $H$. Then
$H$ can be decomposed into an orthogonal sum $H=H_0 \oplus H_1$ of
$V(\cdot)$-invariant subspaces such that the restriction of $V(\cdot)$
on $H_0$ is a unitary (semi)group and the restriction of $V(\cdot)$ on
$H_1$ is a continuous unilateral shift, i.e.,
$H_1$ is unitarily equivalent to $L^2(\R_+, Y)$ for some Hilbert space
$Y$ such that the restriction of $V(\cdot)$ to $H_1$ is equivalent to
the right shift semigroup on $L^2(\R_+, Y)$.
%here exists a subspace
%$Y\subset H_1$ with $V^n Y\perp V^m Y$ for all $n\neq m$, $n,m\in\N$, such that $H_1=\oplus_{n=1}^\infty V^n Y$ holds.
\end{thm}

%As a first application of the Wold decomposition we obtain the density result for periodic operators in  $\mathcal{I}$. (Note that periodic isometries are unitary.)
We further need the following lemma, see also Peller \cite{peller} for the discrete version.

\begin{lemma}\label{lemma:shift-appr}
Let $Y$ be a Hilbert space and let $R(\cdot)$ be the right shift
semigroup on $H:=L^2(\R_+, Y)$. Then there exists a sequence
$\{U_n(\cdot)\}_{n=1}^\infty$ of periodic unitary (semi)groups on $H$
converging strongly to $R(\cdot)$ uniformly on compact time intervals.
\end{lemma}

\begin{proof}
For every $n\in\N$ we define $U_n(\cdot)$ by
\begin{equation*}
(U_n(t)f)(s):=
\begin{cases}
f(s), \quad &s\geq n;\\
R_n(t) f(s), \quad &s\in [0,n],
%f(s-t), \quad &s\in [t,n]; \\
%f(s+n-t), \quad &s\in [0,t].
\end{cases}
\end{equation*}
where $R_n(\cdot)$ denotes the $n$-periodic right shift on the space
$L^2([0,n],Y)$. Then every $U_n(\cdot)$ is a $C_0$-semigroup on $L^2(\R_+,Y)$ which is
isometric and $n$-periodic, and therefore unitary.

Fix $f\in L^2(\R_+,Y)$ 
%and $T>0$. Then for $t\leq T$ and $n>T$ we have
and take $0\leq t<n$. Then we have 
\begin{eqnarray*}
\Vert U_n(t) f - R(t) f\Vert^2
= \int_{n-t}^n \|f(s)\|^2 ds +\int_n^\infty \|f(s)-f(s-t)\|^2 ds %\\
\leq 2 \int_{n-t}^\infty \|f(s)\|^2 ds \underset{n\to\infty}{\longrightarrow} 0
%&=& \int_n^\infty \|f(s)-f(s-t)\|^2 ds + \int_0^t \|f(s+n-t)\|^2 ds \\
%&\leq& \int_n^\infty \|f(s)\|^2 ds + \int_{n-t}^\infty \|f(s)\|^2 ds
%     + \int_{n-t}^n \|f(s)\|^2 ds \\
%&=& 2 \int_{n-t}^\infty \|f(s)\|^2 ds
%   \leq  2 \int_{n-T}^\infty \|f(s)\|^2 ds \underset{n\to\infty}{\longrightarrow} 0
\end{eqnarray*}
uniformly in 
%$t\in [0,T]$ 
$t$ on compact intervals and the lemma is proved.
\end{proof}

As a consequence of the Wold decomposition and Lemma \ref{lemma:shift-appr} we obtain the
following density result for periodic unitary (semi)groups in $\mathcal{I}$.
%and almost weakly stable semigroups in $\mathcal{I}$.

\begin{prop}\label{prop:periodic-i} The set of all periodic unitary semigroups is dense in $\mathcal{I}$.
\end{prop}

\begin{proof} Let $V(\cdot)$ be an isometric semigroup on $H$. Then by
  Theorem \ref{thm:Wold} the orthogonal decomposition $H=H_0\oplus
  H_1$ holds, where the restriction $V_0(\cdot)$ of $V(\cdot)$ to
  $H_0$ is unitary, $H_1$ is unitarily equivalent to
  $L^2(\R_+, Y)$ for some $Y$ and the restriction $V_1(\cdot)$ of
  $V(\cdot)$ on $H_1$ corresponds by this equivalence to the right
  shift semigroup on $L^2(\R_+, Y)$.
By Proposition \ref{prop:periodic} and Lemma \ref{lemma:shift-appr} we can approximate both semigroups $V_0(\cdot)$ and $V_1(\cdot)$ by unitary periodic ones and the assertion follows.
\end{proof}

%As a consequence of the Wold decomposition we obtain the following density result for almost weakly stable operators in $\mathcal{I}$.
We further have a density result for almost weakly stable isometric semigroups. %in $\mathcal{I}$.

\begin{prop}\label{prop:almweakstab-i} The set $\mathcal{W_I}$ of almost weakly stable isometric semigroups is dense in $\mathcal{I}$.
\end{prop}

\begin{proof} Let $V$ be an isometry on $H,\ H_0,\ H_1$ the orthogonal subspaces from Theorem \ref{thm:Wold} and $V_0$ and $V_1$ the corresponding restrictions of $V$. By Lemma \ref{lemma:shift-appr} the operator $V_1$ can be approximated by unitary operators on $H_1$. The assertion now follows from Proposition \ref{prop:almweakstab}.
\end{proof}

By using the same arguments as in the proof of Theorem \ref{thm:unitary} we obtain with the help of Propositions \ref{prop:periodic-i} and \ref{prop:almweakstab-i} the following category result for weakly and almost weakly stable isometries.

\begin{thm} \label{thm:isometry}
The set $\mathcal{S_I}$ of all weakly stable isometric semigroups is
of first category and the set $\mathcal{W_I}$ of all almost weakly
stable isometric semigroups is residual in $\mathcal{I}$.
\end{thm}

%\vspace{0.1cm}

%%%%%%%%%%%%%%%%%%%%%%%%%%%%%%%%%%%%%%%%%%%%%%%%%%%%%%%%%%%%%%%%%%%%%%%
%%                                                                   %%
%%             KONTRAKTIVER  FALL                                    %%
%%                                                                   %%
%%%%%%%%%%%%%%%%%%%%%%%%%%%%%%%%%%%%%%%%%%%%%%%%%%%%%%%%%%%%%%%%%%%%%%%

\section{A remark on the contractive case}

It is not clear how to prove an analogue to Theorems \ref{thm:unitary} and \ref{thm:isometry} for contractive semigroups, while this is done in the discrete case in Eisner, Serény \cite{eisner/sereny:2006}.  We point out one of the difficulties.

Let $\mathcal{C}$ denote the set of all contraction semigroups on $H$ endowed with the metric
\begin{equation*}
d(T(\cdot),S(\cdot)):= \sum_{n,i,j=1}^\infty \frac{\sup_{t\in[0,n]}|\left<T(t)x_i,x_j\right> -\left< S(t)x_i,x_j\right>|}{2^{i+j+n} \|x_i\| \|x_j\|}\quad \text{for }\ T,S \in \mathcal{C},
\end{equation*}
where $\{x_j\}_{j=1}^\infty$ is a fixed dense subset of $H$. The corresponding convergence is the weak convergence of semigroups uniform on compact time intervals. Note that this metric is
a continuous analogue of the metric used in Eisner, Serény \cite{eisner/sereny:2006}.

Since for $H:=l^2$ there exists a Cauchy sequence $\{T_n(\cdot)\}_{n=1}^\infty$ in $\mathcal{C}$ of semigroups (with bounded generators) such that the pointwise limit $S(\cdot)$ does not satisfy the semigroup law (see Eisner, Serény \cite{eisner/sereny:2007a}), the metric space $\mathcal{C}$ is \emph{not} complete (or compact) in general, 
% new
therefore the standard Baire category theorem cannot be applied. It is an open question whether $\mathcal{C}$ is at least \v Cech complete.

\textbf{Acknowledgement.} The authors are very grateful to Rainer Nagel and the referee for valuable comments.

\parindent0pt

\end{document}